\documentclass[twocolumn]{autart}

\usepackage{natbib}
\usepackage{amsmath} 
\usepackage{amssymb}  
\usepackage{mathtools}
\usepackage{mathrsfs}
\usepackage[bbgreekl]{mathbbol}
\usepackage{amsfonts}
\usepackage{color}
\usepackage{psfrag}
\DeclareFontEncoding{LS1}{}{}
\DeclareFontSubstitution{LS1}{stix}{m}{n}
\DeclareSymbolFont{stixsymbols}       {LS1}{stixscr}  {m} {n}
\DeclareSymbolFontAlphabet{\mathscrl} {stixsymbols}

\DeclareSymbolFontAlphabet{\mathbb}{AMSb}
\DeclareSymbolFontAlphabet{\mathbbl}{bbold}
\DeclareMathOperator{\diff}{d}
\providecommand{\ddt}{\tfrac{\diff}{\diff t}}
\newcommand{\kronecker}{\raisebox{1pt}{\ensuremath{\:\otimes\:}}} 
\providecommand{\norm}[1]{\lVert#1\rVert}
\DeclareMathOperator*{\diag}{\text{diag}}

\allowdisplaybreaks

\begin{document}
\begin{frontmatter}

\title{On the steady-state behavior of a nonlinear power system model\thanksref{footnoteinfo}} 

\thanks[footnoteinfo]{This research is supported by ETH funds and the SNF Assistant Professor Energy Grant \#160573. 
A preliminary version of part of the results in in this paper has been presented at the 6th IFAC Work-shop on Distributed Estimation and Control in Networked Systems, September 8-9, 2016, Tokyo, Japan.}

\author[First]{Dominic Gro\ss{}}\ead{gross@control.ee.ethz.ch}
\author[First]{Catalin Arghir}\ead{carghir@control.ee.ethz.ch}
\author[First]{Florian D\"orfler}\ead{dorfler@control.ee.ethz.ch}
 
\address[First]{Automatic Control Laboratory at the Swiss Federal Institute of Technology (ETH) Z\"urich, Switzerland.}

\begin{abstract}                
In this article, we consider a dynamic model of a three-phase power system including nonlinear generator dynamics, transmission line dynamics, and static nonlinear loads. We define a synchronous steady-state behavior which corresponds to the desired nominal operating point of a power system and obtain necessary and sufficient conditions on the control inputs, load model, and transmission network, under which the power system admits this steady-state behavior. We arrive at a separation between the steady-state conditions of the transmission network and generators, which allows us to recover the steady-state of the entire power system solely from a prescribed operating point of the transmission network. Moreover, we constructively obtain necessary and sufficient steady-state conditions based on network balance equations typically encountered in power flow analysis. Our analysis results in several necessary conditions that any power system control strategy needs to satisfy.
\end{abstract}


\end{frontmatter}

\section{Introduction}
The electric power system has been paraphrased as the most complex machine engineered by mankind \citep{K94}. Aside from numerous interacting control loops, power systems are large-scale, and contain highly nonlinear dynamics on multiple time scales from mechanical and electrical domains. As a result power system analysis and control is typically based on simplified models of various degrees of fidelity \citep{PWS-MAP:98}.

A widely accepted reduced power system model is a structure-preserving multi-machine model, where each generator model is reduced to the swing equation modeling the interaction between the generator rotor and the grid, which is itself modeled at quasi-steady-state via the nonlinear algebraic power balance equations, see e.g. \citep{AVDS-TS:16}. Despite being based on  time-scale separations, quasi-stationarity assumptions, and multiple other simplifications, this model has proved itself useful for power system analysis and control \citep{K94,PWS-MAP:98}. Nevertheless, the validity of the simplified model has always been a subject of debate; see \citep{CT15,PM-CDP-NM-AVDS:16} for recent discussions.

The modeling, analysis, and control of power systems has seen a surging research activity in the last years. One particular question of interest concerns the analysis of first-principle nonlinear multi-machine power system models without simplifying generator modeling assumptions and with dynamic (and not quasi-stationary) transmission network models.
\cite{FZO+13} consider a highly detailed power system model based on port-Hamiltonian system modeling, and they carry out a stability analysis for a single generator connected to a constant linear load. This model can be reduced to the classic swing equation model by replacing the electromagnetic generator dynamics with a static relationship between the mechanical power and electrical power supplied to a generator \citep{AVDS-TS:16}.

\cite{CT14} consider stability analysis of a power system using incremental passivity methods. Their analysis requires, among others, the assumptions of a constant torque and field current at the generators. Unfortunately, their analysis also requires a power preservation property that is hard to verify and whose inherent difficulty is rooted in the specific $dq$ coordinates used for the analysis. These coordinates are convenient for a single generator but incompatible for multiple generators \citep{CT:16}.

\cite{BSO+16} study a single generator connected to an infinite bus (i.e. in isolation) and improve upon the previous papers by requiring milder conditions to certify stability, though it is unclear if the analysis can be extended to a multi-machine power systems.
Related stability results have also been obtained for detailed models of so called grid-forming power converters that emulate the dynamics of generators \citep{NW14,JAD17}. Finally, detailed generator models with the transmission networked modeled by quasi-stationary balance equations are studied by \cite{TS-CDP-AVDS:16,DBO+09}. In particular, \cite{DBO+09} study 
existence of equilibria to the nonlinear differential-algebraic model.

In this article we study the port-Hamiltonian power system model derived from first-principles by \cite{FZO+13} and specify the class of steady-state behaviors which are consistent with the nominal operation of a power system, i.e., all three-phase AC signals are balanced, sinusoidal, of constant amplitude, and of the same synchronous frequency. The aim of power systems control is to stabilize such a synchronous steady-state. However, no results are available that give conditions under which detailed first-principles power system models, such as the one proposed by \cite{FZO+13}, admit such a steady-state behavior. Given the importance of the notion of a steady-state for stability analysis and control design, we seek answers to similar questions as in \citep{DBO+09}: under which conditions does the power system admit a synchronous steady-state behavior.

The main contributions of this work are algebraic conditions which relate the state variables, control inputs, load models, and the synchronous steady-state frequency such that the dynamics of the power system coincide with the synchronous steady-state dynamics. We show that the synchronous steady-state behavior is invariant with respect to the power system dynamics if and only if the control inputs and nominal steady-state frequency are constant, as conjectured in \citep{CT14}. We show that load models must be nonlinear ``impedance loads'' to be compatible with the synchronous steady state. Moreover, we obtain a separation between the steady-state conditions of the transmission network and generators, which allows us to explicitly recover the steady state of the entire power system and corresponding inputs from a prescribed steady state of the transmission network. Finally, we constructively obtain conditions based on the well-known network balance (or power flow) equations and show that the power system (with constant inputs) admits a synchronous steady-state behavior if and only if these equations can be solved. 

This paper is organized as follows: In Section \ref{sec:not}, we introduce some basic definitions, the first-principle nonlinear dynamical model of a power system, and the synchronous steady-state dynamics. The results on steady-state operation of the power system are presented in Section \ref{sec:cond}. In Section \ref{sec:disc} we state the main result and discuss its implications. The paper closes with some conclusions in Section \ref{sec:concl}.

\section{Notation and Problem Setup}\label{sec:not}

\subsection{Notation}
We use $\mathbb{R}$ to denote the set of real numbers and $\mathbb{R}_{>0}$ to denote the set of positive real numbers. The set $\mathbb{S}^1$ denotes the unit circle, an angle is a point $\theta \in \mathbb{S}^1$. For column vectors $x \in \mathbb{R}^n$, $y \in \mathbb{R}^m$ we use $(x,y) = [x^{\top} \; y^{\top}]^{\top} \in \mathbb{R}^{n+m}$ to denote a stacked vector, and for vectors or matrices $x$, $y$ we use $\diag(x,y)=\left[\begin{smallmatrix} x & 0 \\ 0 & y\end{smallmatrix}\right]$. Matrices of zeros and ones of dimension $n \times m$ are denoted by $\mathbbl{0}_{n \times m}$ and $\mathbbl{1}_{n \times m}$, and $\mathbbl{0}_n$ and $\mathbbl{1}_n$ denote corresponding column vectors of length $n$. Given $\theta \in \mathbb{S}^1$ we define the rotation matrix $\mathrm{R}(\theta)$, the $90^\circ$ rotation matrix $j$, and the vector $\mathrm{r}(\theta)$ by
\begin{align*}
  \mathrm{R}(\theta) \!\coloneqq\!\! \begin{bmatrix} \cos(\theta)\!\! &-\sin(\theta)\\ \sin(\theta)\!\! & \cos(\theta)\end{bmatrix}\!,\;\;
  j \!\coloneqq\! \mathrm{R}(\pi/2),\;\;
  \mathrm{r}(\theta) \!\coloneqq\!\! \begin{bmatrix} \cos(\theta) \\ \sin(\theta) \end{bmatrix}\!.
\end{align*}
Furthermore, we define the matrix $\mathscrl{j} = \diag(j,\mathbbl{0}_{3 \times 3})$, $I_n$ denotes the identity matrix of dimension $n$, $\otimes$ denotes the Kronecker product, and $\norm{x}=\sqrt{x^{\top}x}$ denotes the Euclidean norm. 

\subsection{Dynamical Model of a Power System}
In this work, we consider a dynamical model of a three-phase power system including nonlinear generator dynamics, transmission line dynamics, and static nonlinear loads derived in \cite{FZO+13}. The reader is referred to \cite{FZO+13} and the references therein for a detailed derivation.
The following assumption is required to prove the main result of the manuscript.
\begin{assum}\label{assum:balanced}
 It is assumed that the three-phase electrical components (resistance, inductance, capacitance) of each device have identical values for each phase. 
\end{assum}
In addition, all three-phase state variables (i.e., voltages and currents) evolve on a two dimensional subspace of the three-phase $abc$-frame, called the $\alpha\beta$-frame, \citep{C43} during synchronous and balanced steady-state operation. Thus, without loss of generality, we can restrict the analysis to the $\alpha\beta$-frame. The dynamical model used throughout this manuscript is obtained by transforming the model published in \cite{FZO+13} into the $\alpha\beta$-frame by applying the Clarke transformation \citep{C43} to each three-phase variable. Moreover, we will work with the {\em co-energy variables}, i.e., voltages and currents, instead of the natural {\em Hamiltonian energy variables}, i.e., charges and fluxes. We emphasize that the change of coordinates is used for simplicity of notation and all of the results also apply to the model derived in \cite{FZO+13} using {\em port-Hamiltonian energy variables} in $abc$-coordinates. 

Because $\alpha\beta$-coordinates can be interpreted as an embedding of the complex numbers into real-valued Euclidean coordinates, the $90^\circ$ rotation matrix $j$ plays the same role that the imaginary unit $\sqrt{-1}$ plays in traditional power system analysis in complex coordinates. We carry out our analysis in the stationary $\alpha\beta$-frame to avoid any limitations which may arise by restricting the analysis to a specific rotating coordinate frame. 

\subsubsection{Power System Topology}
The power system model used in this work consists of $n_g$ generators and $n_v$ AC voltage buses that are interconnected via $n_t$ transmission lines. The topology of the transmission network is described by the oriented incidence matrix $E \in \{-1,1,0\}^{n_v \times n_t}$ of its associated graph, i.e., the $n_v$ voltage buses are the nodes and the $n_t$ transmission lines are the edges of the graph induced by $E$. The incidence matrix of the AC network accounting for $\alpha\beta$-coordinates is denoted by $\mathcal{E} \coloneqq E \otimes I_2 \in \mathbb{R}^{2n_v\times 2n_t}$.

\subsubsection{Transmission Network}
\begin{figure}[t!!!]
\begin{center}
   \psfrag{o}{$v$}
   \psfrag{r}{$v$}   
   \psfrag{p}{$i_T$}
   \psfrag{q}{$L_T$}
   \psfrag{t}{$C$}
   \psfrag{u}{$C$}   
   \psfrag{w}{$R_T$}
   \psfrag{x}{$i_{\textit{in}}$}   
   \includegraphics[trim=0 0cm 0 0, clip, width=5cm]{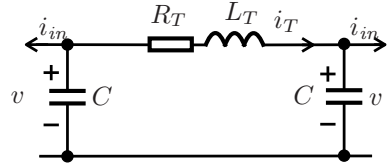}
   \caption{Transmission line connecting two voltage buses.\label{fig:trans}}
   \end{center}
\end{figure}

The transmission lines are modeled using the $\Pi$-model \citep{PWS-MAP:98} depicted in Figure \ref{fig:trans}. The state variables of the power grid model are the $n_t$ line currents $i_T \coloneqq (i_{T,1},\ldots,i_{T,n_t}) \in \mathbb{R}^{2n_t}$ and $n_v$ bus voltages $v \coloneqq (v_1,\ldots,v_{n_v}) \in \mathbb{R}^{2n_v}$, where $i_{T,\mathsf{k}} \in \mathbb{R}^2$ is the current flowing through a transmission line $\mathsf{k} \in \{1,\ldots,n_t\}$ and $v_\mathsf{k} \in \mathbb{R}^2$ is the voltage of a voltage bus $\mathsf{k} \in \{1,\ldots,n_v\}$. The transmission network is interconnected to the loads and generators via the currents $i_{\textit{in}} \coloneqq (i_{\textit{in},1},\ldots,i_{\textit{in},n_v})  \in \mathbb{R}^{2 n_v}$, where $i_{\textit{in},\mathsf{k}} \in \mathbb{R}^2$ is the current flowing out of an AC voltage bus $\mathsf{k} \in \{1,\ldots,n_v\}$ and into a generator or a load. The dynamics of the $n_v$ voltage buses and $n_t$ transmission lines are given by
\begin{subequations}\label{eq:network}
\begin{align}
C \ddt v &= - \mathcal{E} i_T - i_{\textit{in}}, \label{eq:vbus} \\
L_T \ddt i_T &= -R_T i_T + \mathcal{E}^{\top} v, \label{eq:trline}
\end{align} 
\end{subequations}
where $C \coloneqq \diag(C_1,\ldots,C_{n_v}) \in \mathbb{R}^{2n_v \times 2 n_v}$ is the matrix of voltage bus capacitances $C_\mathsf{k}=I_2 c_\mathsf{k}$ with $c_\mathsf{k} \in \mathbb{R}_{>0}$, $L_T \coloneqq \diag(L_{T,1},\ldots,L_{T,n_t}) \in \mathbb{R}^{2n_t \times 2n_t}$ aggregates the line inductances $L_{T,\mathsf{k}}=I_2 l_{T,\mathsf{k}}$ with $l_{T,\mathsf{k}} \in \mathbb{R}_{>0}$, and $R_T \coloneqq \diag(R_{T,1},\ldots,R_{T,n_t}) \in \mathbb{R}^{2n_t \times 2n_t}$ is the matrix of line resistances $R_{T,\mathsf{k}}=I_2 r_{T,\mathsf{k}}$ with $r_{T,\mathsf{k}} \in \mathbb{R}_{>0}$. 

\subsubsection{Synchronous Machines}
\begin{figure}[t!!!]
\begin{center}
   \psfrag{g}{$\tau_m$}
   \psfrag{h}{$\theta$, $\omega$}
   \psfrag{i}{$v_f$}
   \psfrag{x}{$\theta$}   
   \psfrag{k}{$i_f$}
   \psfrag{d}{$i_d$}
   \psfrag{q}{$i_q$}
   \psfrag{e}{$r_q$}
   \psfrag{p}{$r_d$}
   \psfrag{j}{$v_{\alpha}$}
   \psfrag{jj}{$v_{\beta}$}
   \psfrag{l}{$\tau_e$}
   \psfrag{m}{$i_{\alpha}$}
   \psfrag{u}{$i_{\beta}$}
   \psfrag{w}{$L(\theta)$}
   \psfrag{n}{}
   \psfrag{z}{$M$}
   \psfrag{y}{$r_f$}
   \psfrag{rs}{$r_s$}   
   \includegraphics[trim=0 0cm 0 0, clip, width=0.95\columnwidth]{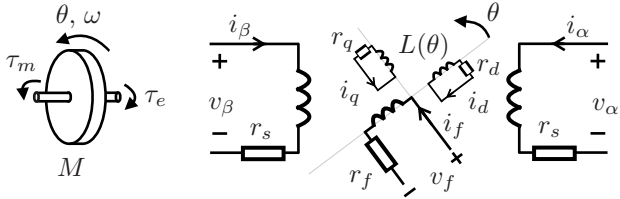}
   \caption{Mechanical and electrical components of a synchronous machine.\label{fig:syncmachine}}
   \end{center}
\end{figure}

A synchronous machine with index $\mathsf{k} \in \{1,\ldots,n_g\}$ is modeled by
\begin{subequations}%
 \label{eq:gendyn}
\begin{align}
\ddt \theta_\mathsf{k} &= \omega_\mathsf{k} \\
m_\mathsf{k} \ddt \omega_\mathsf{k} &= -d_\mathsf{k} \omega_\mathsf{k} - \tau_{e,\mathsf{k}} + \tau_{m,\mathsf{k}}\\
L_\mathsf{k}(\theta_\mathsf{k}) \ddt i_\mathsf{k} &= -R_\mathsf{k} i_\mathsf{k} +  \left[\begin{smallmatrix}  v_\mathsf{k} \\ v_{f,\mathsf{k}} \\ \mathbbl{0}_2 \end{smallmatrix}\right] - v_{\textit{ind},\mathsf{k}}, \label{eq:gendyn:elec}
\end{align}
\end{subequations}
where $i_\mathsf{k} \!=\! (i_{s,\mathsf{k}},i_{r,\mathsf{k}}) \!\in\! \mathbb{R}^5$ aggregates the stator currents $i_{s,\mathsf{k}}\!=\!(i_{\alpha,\mathsf{k}},i_{\beta,\mathsf{k}}) \!\in\! \mathbb{R}^2$ and rotor currents $i_{r,\mathsf{k}}\!=\!(i_{f,\mathsf{k}},i_{d,\mathsf{k}},i_{q,\mathsf{k}}) \!\in\! \mathbb{R}^3$, with excitation current $i_{f,\mathsf{k}}$, and the damper winding currents $i_{d,\mathsf{k}}$ and $i_{q,\mathsf{k}}$. Moreover, $v_\mathsf{k} = (v_{\alpha,\mathsf{k}}, v_{\beta,\mathsf{k}}) \in \mathbb{R}^2$ denotes the voltages at the generators voltage bus, $\omega_\mathsf{k} \in \mathbb{R}$ is rotational speed of the rotor, $\theta_\mathsf{k} \in \mathbb{S}^1$ its angular displacement, $\tau_{e,\mathsf{k}} \in \mathbb{R}$ is the electrical torque acting on the rotor, and $v_{\textit{ind},\mathsf{k}} \in \mathbb{R}^5$ is a voltage induced by the rotation of the machine. The machine is actuated by the voltage $v_{f,\mathsf{k}} \in \mathbb{R}$ across the excitation winding of the generator and the mechanical torque $\tau_{m,\mathsf{k}} \in \mathbb{R}$ applied to the rotor. The mechanical and electrical part of the machine are depicted in Figure \ref{fig:syncmachine}. The inertia and damping of the rotor are denoted by $m_\mathsf{k} \!\in \mathbb{R}_{>0}$ and $d_\mathsf{k} \!\in \mathbb{R}_{>0}$, and the windings have resistance $R_\mathsf{k}\! =\! \diag(R_{s,\mathsf{k}},R_{r,\mathsf{k}})$ with stator resistance $R_{s,\mathsf{k}}\!=I_2 r_{s,\mathsf{k}}$, $r_{s,\mathsf{k}}\!\in \mathbb{R}_{>0}$ and rotor resistance $R_{r,\mathsf{k}}\!=\!\diag(r_{f,\mathsf{k}},r_{d,\mathsf{k}},r_{q,\mathsf{k}})$, with excitation winding resistance $r_{f,\mathsf{k}} \in \mathbb{R}_{>0}$, and damper winding resistances $r_{d,\mathsf{k}} \in \mathbb{R}_{>0}$ and $r_{q,\mathsf{k}} \in \mathbb{R}_{>0}$. The inductance matrix $L_\mathsf{k}(\theta_\mathsf{k}): \mathbb{S}^1 \to \mathbb{R}^{5 \times 5}$ is defined using the stator inductance $L_{s,\mathsf{k}}(\theta_\mathsf{k})\in \mathbb{R}^{2 \times 2}$, rotor inductance $L_{r,\mathsf{k}} \in \mathbb{R}^{3 \times 3}$, and mutual inductance $L_{m,\mathsf{k}}(\theta_\mathsf{k})\in \mathbb{R}^{2 \times 3}$:
\begin{align*}
L_\mathsf{k}(\theta_\mathsf{k}) \!=\!\! \begin{bmatrix} L_{s,\mathsf{k}}(\theta_\mathsf{k}) & L_{m,\mathsf{k}}(\theta_\mathsf{k})\! \\  
                                    \!L_{m,\mathsf{k}}(\theta_\mathsf{k})^{\!\top} &  L_{r,\mathsf{k}} \end{bmatrix}\!, \;
L_{r,\mathsf{k}} \!=\!\! \begin{bmatrix} l_{f,\mathsf{k}} & l_{fd,\mathsf{k}} & 0 \\ l_{fd,\mathsf{k}} & l_{d,\mathsf{k}} & 0 \\ 0 & 0 & l_{q,\mathsf{k}} \end{bmatrix}\!\!,
\end{align*}
where $l_{f,\mathsf{k}} \in \mathbb{R}_{>0}$, $l_{d,\mathsf{k}} \in \mathbb{R}_{>0}$, $l_{q,\mathsf{k}} \in \mathbb{R}_{>0}$, and $l_{fd,\mathsf{k}} \in \mathbb{R}_{>0}$ are the inductances of the excitation winding, damper windings, and the mutual inductance of the excitation and damper winding, respectively. With the stator winding inductance $l_{s,\mathsf{k}} \in \mathbb{R}_{>0}$ and rotor saliency $l_{sa,\mathsf{k}}\in \mathbb{R}_{\geq0}$ the stator inductance is given by
\begin{align*}
 L_{s,\mathsf{k}}(\theta_\mathsf{k}) =  \begin{bmatrix} l_{s,\mathsf{k}} & 0 \\ 0 & l_{s,\mathsf{k}} \end{bmatrix} +  \mathrm{R}(2\theta_\mathsf{k}) \begin{bmatrix} l_{sa,\mathsf{k}} & 0 \\ 0 & -l_{sa,\mathsf{k}}\end{bmatrix}.
\end{align*}
Finally, $L_{m,\mathsf{k}}(\theta_\mathsf{k})$ is defined using the mutual inductances $l_{sf,\mathsf{k}} \in \mathbb{R}_{>0}$, $l_{sd,\mathsf{k}} \in \mathbb{R}_{>0}$, $l_{sq,\mathsf{k}} \in \mathbb{R}_{>0}$:
\begin{align*}
 L_{m,\mathsf{k}}(\theta_\mathsf{k}) =  \mathrm{R}(\theta_\mathsf{k}) \begin{bmatrix} l_{sf,\mathsf{k}} & l_{sd,\mathsf{k}} & 0 \\ 0 & 0 &-l_{sq,\mathsf{k}} \end{bmatrix}.
\end{align*} 
The electrical torque acting on the rotor is given by
\begin{align}
 \tau_{e,\mathsf{k}} = \tfrac{1}{2} i^{\top}_\mathsf{k} \big(L_\mathsf{k}(\theta_\mathsf{k}) \mathscrl{j}  + \mathscrl{j}^{\top} L_\mathsf{k}(\theta_\mathsf{k})\big) i_\mathsf{k}, \label{orig:torque}
\end{align}
and the voltage $v_{\textit{ind},\mathsf{k}} \in \mathbb{R}^5$ induced in the machine windings due to the rotation of the machine is given by
\begin{align}\label{eq:vind}
 v_{\textit{ind},\mathsf{k}} = \omega_\mathsf{k} \big(L_\mathsf{k}(\theta_\mathsf{k}) \mathscrl{j}^\top + \mathscrl{j} L_\mathsf{k}(\theta)\big) i_\mathsf{k}.
\end{align}

\subsubsection{Static Loads}\label{sec:loads}
In this work, we consider the static load model used in \cite{FZO+13}. Specifically, loads are included in the model via a load current $i_{l,\mathsf{k}}: \mathbb{R}^2 \to \mathbb{R}^2$ (flowing out of a voltage bus $\mathsf{k} \in \{1,\ldots,n_v\}$) that is a function of the bus voltage $v_\mathsf{k}$, and satisfies the dissipation inequality $i_{l,\mathsf{k}}(v_\mathsf{k})^\top v_\mathsf{k} \geq 0$. We additionally assume that $i_{l,\mathsf{k}}(v_\mathsf{k}) = \mathbbl{0}_2$ if $v_\mathsf{k}=\mathbbl{0}_2$.

\subsubsection{Dynamic Model of the Power System}
To obtain the dynamic model of the entire power system, the synchronous machine model \eqref{eq:gendyn} and load model are combined with the transmission network model \eqref{eq:network} by defining $i_{\textit{in}} \coloneqq ( i_s, \mathbbl{0}_{2n_l}) + i_l$, where $n_l = n_v - n_g$, $i_l \coloneqq \big(i_{l,1}(v_1),\ldots,i_{l,n_v}(v_{n_v})\big) \in \mathbb{R}^{2n_v}$, and $i_s \coloneqq \big(i_{s,1},\ldots,i_{s,n_g}\big) \in \mathbb{R}^{2n_g}$. With the aggregated vectors $\theta \coloneqq (\theta_1, \ldots, \theta_{n_g})$, $ \omega \coloneqq ( \omega_1,  \ldots, \omega_{n_g})$,  $i \coloneqq (i_1, \ldots,  i_{n_g})$, the state vector of the nonlinear power system model is $x=(\theta, \omega, i, v, i_T) \in \mathbb{R}^{n_x}$, with $n_x = 7n_g+2n_v+2n_T$. Using the vectors $\tau_m \coloneqq (\tau_{m,1}, \ldots, \tau_{m,n_g})$, and $v_f   = (v_{f,1}, \ldots, v_{f,n_g})$ the inputs are given by $u=(\tau_m, v_f) \in \mathbb{R}^{n_u}$, $n_u = 2n_g$. Moreover, to simplify the notation, we define $\tau_e \coloneqq (\tau_{e,1},\ldots, \tau_{e,n_g})$, $v_{\textit{ind}} \coloneqq (v_{\textit{ind},1},\ldots, v_{\textit{ind},n_g})$, the matrices $M$, $D$, $R$, $R_T$, $L(\theta)$, and $L_T$, which collect the node matrices (e.g., $M = \diag(M_1,\ldots,M_{n_g})$), and $\mathcal{M}(x)=\diag(I_{n_g},M,L(\theta),C,L_T)$ collecting the time constants. We use indicator matrices $\mathcal{I}_{f} = I_{n_g} \kronecker (0, 0,  1, 0, 0)$, $\mathcal{I}_s = I_{n_g} \kronecker (I_2,  \mathbbl{0}_{3 \times 2})$, and
$\mathcal{I}^{\top}_{v} =[ \mathcal{I}_s \; \mathbbl{0}_{5 n_g \times 2n_l}]$ to describe the interconnection of the components results in the following model of the entire power system:
\begin{align}\label{eq:nlmodel}
 \ddt {x} = \mathcal{M}(x)^{-1} \begin{bmatrix}
 \omega \\
-D \omega - \tau_e + \tau_m\\
-R i + \mathcal{I}^{\top}_v v + \mathcal{I}_{f} v_f - v_{\textit{ind}}\\
- \mathcal{I}_v i - \mathcal{E} i_T - i_l \\
-R_T i_T + \mathcal{E}^{\top} v
 \end{bmatrix} = f(x,u),\raisetag{12pt}
\end{align}

\subsection{Synchronous Steady-State Behavior}\label{subsec:desssdyn}
We formulate the following steady-state dynamics to describe synchronous and balanced steady-state operation of a power system at a constant frequency $\omega_0 \in \mathbb{R}$ as
\begin{subequations}%
\label{eq:res}%
\begin{align}%
 \ddt \theta_\mathsf{k} &= \omega_0, &\forall \mathsf{k} \in \{1,\ldots,n_g\}, \label{res:ode:theta}\\
 \ddt \omega_\mathsf{k} &= 0, &\forall \mathsf{k} \in \{1,\ldots,n_g\}, \label{res:ode:p}\\
 \ddt i_\mathsf{k} &=\omega_0\mathscrl{j} i_\mathsf{k}, &\forall \mathsf{k} \in \{1,\ldots,n_g\}, \label{res:ode:lambda}\\ 
 \ddt v_\mathsf{k} &=\omega_0j v_\mathsf{k}, &\forall \mathsf{k} \in \{1,\ldots,n_v\}, \label{res:ode:q}\\
 \ddt i_{T,\mathsf{k}} &=\omega_0j i_{T,\mathsf{k}}, &\forall \mathsf{k} \in \{1,\ldots,n_t\}. \label{res:ode:lambdag}
\end{align}%
\end{subequations}%
The steady-state behavior \eqref{eq:res} specifies a balanced, synchronous, and sinusoidal operation of each grid component.
This results in the steady-state dynamics $\ddt {x} = f_d(x,\omega_0)$ with constant nominal frequency $\omega_0 \in \mathbb{R}$ and
\begin{align}
\label{eq:ssmodel}
  f_d(x,\omega_0) \!=\!      ( \mathbbl{1}_{n_g} \omega_0,\, \mathbbl{0}_{n_g},\,\omega_0\mathcal{J}_g i,\,\omega_0J_{v} v,\,\omega_0J_T i_T),
\end{align}
where $J_v = I_{n_v} \otimes j$, $J_T = I_{n_t} \otimes j$, and $\mathcal{J}_g = I_{n_g} \otimes \mathscrl{j}$. 

In the next section, we derive necessary and sufficient conditions under which the nonlinear power system dynamics \eqref{eq:nlmodel} coincide with the steady-state dynamics \eqref{eq:ssmodel} for all time. In particular, we derive conditions on the control inputs and load models, and obtain a separation between the steady-state conditions of the transmission network and generators. This result allows to explicitly recover a steady state of the entire power system and corresponding steady-state inputs from a prescribed steady state of the transmission network. As we will see, this result also implies that the power system admits a non-trivial synchronous steady-state dynamics \eqref{eq:ssmodel} if and only if there exists a non-trivial solution to the well-known power flow equations.

\section{Conditions for the Existence of Synchronous Steady-States}\label{sec:cond}
We begin our analysis by defining the set $\mathcal{S}$ on which the vector fields of the dynamics \eqref{eq:nlmodel} and \eqref{eq:ssmodel} coincide. In other words, the residual dynamics $\rho(x,u,\omega_0) \coloneqq \mathcal{M}(x) \big( f_d(x,\omega_0) - f(x,u) \big)$ vanish on $\mathcal{S}$:
\begin{align*}
\mathcal{S} \coloneqq \left\{ (x,u,\omega_0) \in \mathbb{R}^{n_x} \times \mathbb{R}^{n_u} \times \mathbb{R} \;\left\vert\; \rho(x,u,\omega_0)=\mathbbl{0}_{n_x} \right. \right\}.
\end{align*}
In the remainder we will give conditions under which the steady-state exists, i.e., under which $\mathcal{S}$ is non-empty, Moreover, we characterize the inputs $u$ and load currents $i_l$ such that $\mathcal{S}$ is invariant, i.e., that trajectories of \eqref{eq:nlmodel} starting in $\mathcal{S}$ coincide with the synchronous steady-state behavior \eqref{eq:ssmodel} for all times.

In the next subsection, we will first derive necessary and sufficient conditions on $u(t)$ and the load model under which $\mathcal{S}$ is invariant, i.e., all trajectories of the dynamics \eqref{eq:nlmodel} starting in $\mathcal{S}$ remain in $\mathcal{S}$ for all time. We thereby characterize the steady-state control inputs and class of static load models, for which \eqref{eq:ssmodel} is a steady-state behavior of the power system \eqref{eq:nlmodel}. Based on these results we will provide conditions for the existence of states $x$ and inputs $u$ such that $(x,u,\omega_0) \in \mathcal{S}$. Together, these conditions are necessary and sufficient conditions for the power system to admit the synchronous steady-state behavior \eqref{eq:ssmodel}.

\subsection{Invariance of the Set $\mathcal{S}$}
To establish invariance of $\mathcal{S}$, we consider the dynamics obtained by combining the nonlinear power grid dynamics described by $f: \mathbb{R}^{n_x} \times \mathbb{R}^{n_u} \to \mathbb{R}^{n_x}$ with controller (torque and excitation) dynamics $g_u: \mathbb{R}_{>0} \times \mathbb{R}^{n_x} \times \mathbb{R}^{n_u} \times \mathbb{R} \to \mathbb{R}^{n_u}$ possibly both depending on time, the system state, and the nominal synchronous frequency $\omega_0$:
\begin{align}\label{eq:exDyn}
\ddt \big(x,u\big) = \big(f(x,u), g_u(t,x,u,\omega_0)\big).
\end{align}

Moreover, the residual dynamics $\rho(x,u,\omega_0)$ can be rewritten as follows:
\begin{align*}
\rho(x,u,\omega_0) = \begin{bmatrix}
   \mathbbl{1}_{n_g} \omega_0 - \omega\\
   D \omega + \tau_e - \tau_m\\
     \big(R + \omega_0 L(\theta) \mathcal{J}_g \big) i - \mathcal{I}^{\top}_v v  - \mathcal{I}_f  v_f + v_{\textit{ind}}\\
 \omega_0C  J_{v} v + \mathcal{I}_v i + \mathcal{E} i_T + i_l \\
  \big(R_T  +\omega_0L_T J_T \big) i_T - \mathcal{E}^{\top} v
 \end{bmatrix}.
\end{align*}

\begin{thm}{\bf(Invariance condition)}\label{th:invariance}
The set $\mathcal{S}$ is invariant with respect to the dynamics \eqref{eq:exDyn} if and only if the control inputs are constant on $\mathcal{S}$, i.e., $g_u(t,x,u,\omega_0)=\mathbbl{0}_{n_u}$ for all $(x,u,\omega_0) \in \mathcal{S}$, and the load current $i_l(v)$ satisfies $\ddt i_l(v) =\omega_0J_v i_l(v)$ for all $(x,u,\omega_0) \in \mathcal{S}$.
\end{thm}
\begin{pf} 
The set $\mathcal{S}$ is invariant with respect to \eqref{eq:exDyn} if and only if $\ddt \rho(x,u,\omega_0) = \mathbbl{0}_{n_x}$ for all $(x,u,\omega_0) \in \mathcal{S}$. On $\mathcal{S}$ it holds that $\ddt x = f_d(x,\omega_0)$, resulting in the following necessary and sufficient condition for invariance of $\mathcal{S}$ for all $(x,u,\omega_0) \in \mathcal{S}$:
 \begin{align}\label{eq:invcond}
 \frac{\diff \rho}{\diff t} =  \frac{\partial \rho}{\partial x} f_d + \frac{\partial \rho}{\partial u} g_u = \mathbbl{0}_{n_x}.
 \end{align}
In the remainder of the proof we derive conditions on the control inputs $u$ and load current $i_l$ that ensure that \eqref{eq:invcond} holds. Using $\mathscrl{j}=-\mathscrl{j}^{\top}$ it can be verified that
\begin{align*}
 \frac{\partial \tau_{e,\mathsf{k}}}{\partial \theta_\mathsf{k}}\omega_0+ \frac{\partial \tau_{e,\mathsf{k}}}{\partial i_\mathsf{k}}\omega_0\mathscrl{j} i_\mathsf{k} &= 0, \\
 \frac{\partial v_{\textit{ind},\mathsf{k}}}{\partial \theta_\mathsf{k}}\omega_0+ \frac{\partial v_{\textit{ind},\mathsf{k}}}{\partial i_\mathsf{k}}\omega_0\mathscrl{j} i_\mathsf{k} &=\omega_0\mathscrl{j} v_{\textit{ind},\mathsf{k}},
\end{align*}
and using $\left(\frac{\partial L(\theta)\mathcal{J}_g i}{\partial x} \right) f_d =\omega_0\mathcal{J}_g L(\theta) \mathcal{J}_g i$ one obtains
\begin{align*}
\! \left(\frac{\partial \rho}{\partial x}\right)\! f_d \! =  \!\omega_0\! \begin{bmatrix}
\mathbbl{0}_{2n_g}\\
\!(R \mathcal{J}_g +\omega_0\mathcal{J}_g L(\theta) \mathcal{J}_g ) i - \mathcal{I}^{\top}_v J_v v + \mathcal{J}_g v_{\textit{ind}}\\
\omega_0 C J^2_v v + \mathcal{I}_v  \mathcal{J}_g i + \mathcal{E} J_T i_T + \tfrac{\partial i_l}{\partial v} J_v v \\
(R_T +\omega_0L_T J_T ) J_T i_T - \mathcal{E}^{\top} J_v v
\end{bmatrix}\!\!.
\end{align*}
Next, it can be verified that $\mathcal{I}^{\top}_v J_v = \mathcal{J}_g \mathcal{I}^{\top}_v$, and $\mathcal{E} J_T = J_{v} \mathcal{E}$. Moreover, $J_v$, $J_T$ and $\mathcal{J}_g$ commute with diagonal matrices. This results in
\begin{align*}
\left(\frac{\partial \rho}{\partial x}\right)\! f_d \! = \omega_0\! \begin{bmatrix}
\mathbbl{0}_{2n_g}\\
\mathcal{J}_g \big((R  +\omega_0L(\theta) \mathcal{J}_g ) i - \mathcal{I}^{\top}_v v + v_{\textit{ind}} \big)\\
J_v \big(\omega_0C J_v v + \mathcal{I}_v i + \mathcal{E} i_T\big) + \tfrac{\partial i_l}{\partial v} J_v v \\
J_T \big((R_T +\omega_0L_T J_T )  i_T - \mathcal{E}^{\top} v\big)
\end{bmatrix}\!.
\end{align*}
By definition of $\mathcal{S}$ we have for all $(x,u,\omega_0) \in \mathcal{S}$ that
\begin{align}\label{eq:drhodx}
&\!\!\!\!\!\left(\frac{\partial \rho}{\partial x}\right) f_d \!=\!  \omega_0\! \begin{bmatrix}
\mathbbl{0}_{2n_g}\\
\mathcal{J}_g \mathcal{I}_f v_f \\
\tfrac{\partial i_l}{\partial v} J_v v \!-\! J_v i_l \\
\mathbbl{0}_{2n_t}
\end{bmatrix} 
\! = \!
\omega_0 \! \begin{bmatrix}
\mathbbl{0}_{2n_g}\\
\mathbbl{0}_{2n_g} \\
\tfrac{\partial i_l}{\partial v} J_v v \!-\! J_v i_l \\
\mathbbl{0}_{2n_t}
\end{bmatrix}\!,\!
\end{align}
where the last equality follows from $\mathcal{J}_g \mathcal{I}_f=0$. 
Using \eqref{eq:invcond} and \eqref{eq:drhodx} one obtains the following condition for invariance of $\mathcal{S}$ with respect to the dynamics \eqref{eq:exDyn} for all $(x,u,\omega_0) \in \mathcal{S}$:
\begin{align*}
 \begin{bmatrix} \mathbbl{0}_{n_g \times n_g} & \mathbbl{0}_{n_g \times n_g} \\ 
  I_{n_g} & \mathbbl{0}_{n_g \times n_g} \\ 
  \mathbbl{0}_{3n_g \times n_g} & \mathcal{I}_f \\ 
  \mathbbl{0}_{2n_v \times n_g} & \mathbbl{0}_{2n_v \times n_g} \\
\mathbbl{0}_{2n_T \times n_g} & \mathbbl{0}_{2n_T \times n_g} 
 \end{bmatrix}\!
  g_u \! = \omega_0\! \begin{bmatrix}
\mathbbl{0}_{n_g}\\
\mathbbl{0}_{n_g}\\
\mathbbl{0}_{2n_g} \\
J_v i_l  \!-\! \tfrac{\partial i_l}{\partial v} J_v v\\
\mathbbl{0}_{2n_t}
\end{bmatrix}\!.
\end{align*}
This holds if and only if $g_u=\mathbbl{0}_{n_u}$, and $\tfrac{\partial i_l}{\partial v} \omega_0 J_v v = \omega_0 J_v i_l$. Moreover, on $\mathcal{S}$ it holds that $\ddt v =\omega_0J_v v$, and it follows that $\tfrac{\partial i_l}{\partial v} \omega_0 J_v v = \tfrac{\partial i_l}{\partial v} \ddt v = \ddt i_l$, i.e., on $\mathcal{S}$ the conditions $\tfrac{\partial i_l}{\partial v} \omega_0 J_v v = \omega_0 J_v i_l$ and $\ddt i_l =\omega_0 J_v i_l$ are identical. \hfill $\square$
\end{pf}
Theorem \ref{th:invariance} shows that operation at the synchronous steady-state behavior \eqref{eq:ssmodel} requires constant torque $\tau_m$ and constant excitation voltage $v_f$ on $\mathcal{S}$. Moreover, on $\mathcal{S}$, the load current $i_{l,\mathsf{k}}(v_\mathsf{k})$ needs to be synchronous and of constant amplitude.
\begin{rem}{\bf(Constant synchronous frequency)}
Using the same approach as in the proof of Theorem \ref{th:invariance} it is straightforward to show that the synchronous frequency $\omega_0$ necessarily has to be constant.
\end{rem}

Based on this result we restrict ourselves from now to a constant (possibly zero) frequency $\omega_0 \in \mathbb{R}$ corresponding to the nominal synchronous operating frequency of the 
power system and constant inputs, i.e., $g_u \coloneqq \mathbbl{0}_{n_u}$. This results in the following set parametrized in $\omega_0 \in \mathbb{R}$:
\begin{align}\label{eq:set:Sw0}
 \!\!\mathcal{S}_{\omega_0} \!:=\! \left\{  (x,u) \in \mathbb{R}^{n_x} \times \mathbb{R}^{n_u} \;\left\vert\; \rho(x,u,\omega_0) = \mathbbl{0}_{n_x}  \right.\right\}\!.\!
\end{align}
In the next subsection, we will use Theorem \ref{th:invariance} to characterize the class of static load models that are consistent with the synchronous steady-state behavior \eqref{eq:ssmodel}.

\subsection{Steady-State Conditions and Load Model}\label{sec:load}
The following result identifies the class of all load models for which the set $\mathcal{S}$ is invariant. In particular, it shows that all load models that are consistent with the synchronous steady-state behavior \eqref{eq:ssmodel} can be expressed as nonlinear impedance loads whose impedance is a function of the voltage magnitude.
\begin{thm}{\bf(Consistent load model)}\label{th:load}
Any static load model, i.e., any $i_{l,\mathsf{k}}(v_\mathsf{k})$ that satisfies $i_{l,\mathsf{k}}(\mathbbl{0}_2)= \mathbbl{0}_2$ and $i_{l,\mathsf{k}}(v_\mathsf{k})^\top v_\mathsf{k} \geq 0$, which is consistent with the synchronous steady state, i.e., satisfies the steady-state conditions of Theorem \ref{th:invariance}, can be expressed in the form $i_{l,\mathsf{k}} = Y_{l,\mathsf{k}}(\norm{v_\mathsf{k}}) v_\mathsf{k}$, where $Y_{l,\mathsf{k}} \coloneqq I_2 g_\mathsf{k}(\norm{v_\mathsf{k}}) + j b_\mathsf{k}(\norm{v_\mathsf{k}})$, $g_\mathsf{k}: \mathbb{R}_{\geq 0} \to \mathbb{R}_{\geq 0}$, and $b_\mathsf{k}: \mathbb{R}_{\geq 0} \to \mathbb{R}$. 
\end{thm}
\begin{pf}
Without loss of generality any $i_l(v)$ that satisfies $i_l(v)=\mathbbl{0}_2$ if $v=\mathbbl{0}_2$ can be expressed as $i_{l,\mathsf{k}} = Y_{l,\mathsf{k}}(v_\mathsf{k}) v_\mathsf{k}$. According to Theorem \ref{th:invariance}, consistency of the load model with the steady state requires that $\ddt i_{l,\mathsf{k}}(v_\mathsf{k}) = \omega_0 j i_{l,\mathsf{k}}(v_\mathsf{k})$ holds on $\mathcal{S}$. Moreover, on $\mathcal{S}$ it holds that $\ddt v_\mathsf{k} = \omega_0 j v_\mathsf{k}$. Therefore, in steady-state, the load current $i_{l,\mathsf{k}}(v_\mathsf{k})$ and voltage $v_\mathsf{k}$ are required to have constant amplitude and rotate with the synchronous frequency $\omega_0$, i.e. the relative angle between $i_{l,\mathsf{k}}(v_\mathsf{k})$ and $v_\mathsf{k}$ and the ratio of their amplitudes is required to be constant. It follows that, for every voltage magnitude $\norm{v_\mathsf{k}}$, all load currents $i_{l,\mathsf{k}}$ consistent with the synchronous steady state can be expressed using a relative angle $\delta_\mathsf{k}: \mathbb{R}_{\geq 0} \to \mathbb{S}^1$ and a gain $\mu_\mathsf{k}: \mathbb{R}_{\geq 0} \to \mathbb{R}$ as follows:
\begin{align}\label{eq:loadpol}
 i_{l,\mathsf{k}}(v_\mathsf{k}) = \mu(\norm{v_\mathsf{k}}) \mathrm{R}\big(\delta_\mathsf{k}(\norm{v_\mathsf{k}})\big) v_\mathsf{k}.
\end{align}
Note, that $\mu_\mathsf{k}$ and $\delta_\mathsf{k}$ are required to be constant on $\mathcal{S}$. Because $v_\mathsf{k}$ is not constant on $\mathcal{S}$, the gain $\mu_\mathsf{k}$ and relative angle $\delta_\mathsf{k}$ can only depend on $\norm{v_\mathsf{k}}$. Because $\norm{v_\mathsf{k}}$ is constant on $\mathcal{S}$, it follows that $\mu(\norm{v_\mathsf{k}})$ and $\delta_\mathsf{k}(\norm{v_\mathsf{k}})$ are constant on $\mathcal{S}$. Moreover, for any given angle $\delta_\mathsf{k}$ and gain $\mu_\mathsf{k}$ the polar parametrization of the load model \eqref{eq:loadpol} is equivalent to the load model given in the theorem with $b_\mathsf{k} = \mu_\mathsf{k} \sin(\delta_\mathsf{k})$ and $g_\mathsf{k} = \mu \cos(\delta_\mathsf{k})$. The theorem follows by noting that the dissipation inequality $i_{l,\mathsf{k}}(v_\mathsf{k})^\top v_\mathsf{k} \geq 0$ holds for this load model if and only if $g_\mathsf{k}(\norm{v_\mathsf{k}}) \geq 0$. \hfill $\square$
\end{pf}
The load model identified in Theorem \ref{th:load} allows us to define all static load models commonly used in power system analysis, e.g., constant current loads, constant power loads, and constant impedance loads. In particular, $P_\mathsf{k} = g_\mathsf{k}(\norm{v_\mathsf{k}}) v^\top_\mathsf{k} v_\mathsf{k} $ and $Q_\mathsf{k} = -b_\mathsf{k}(\norm{v_\mathsf{k}}) v^\top_\mathsf{k} v_\mathsf{k}$ are the active and reactive power drawn by the load. Based on this result we can restrict the analysis to static load models according to Theorem \ref{th:load}. To this end, we define $Y_l(v) = \diag(Y_{l,1}(\norm{v_1}),\ldots,Y_{l,n_v}(\norm{v_{n_v}})$ and further specify the load current as $i_l(v) \coloneqq Y_l(v) v$.

Theorem \ref{th:invariance} and Theorem \ref{th:load} show that besides $(x_0,u_0) \in \mathcal{S}_{\omega_0}$ synchronous steady-state operation also requires that the steady-state control inputs $u$ are constant and identifies the class of static load models which are consistent with the synchronous steady state. Next, we derive conditions for the existence of states $x$ and inputs $u$ such that $(x,u) \in \mathcal{S}_{\omega_0}$. 

\subsection{Steady-State Analysis and Network Equations}
In the following, we establish a connection between the steady-state conditions \eqref{eq:set:Sw0} and Kirchhoff's equations for the transmission network (also often formulated as power balance/flow equations). The power system \eqref{eq:nlmodel} can be partitioned into $n_g$ generators with dynamics \eqref{eq:gendyn} and the network consisting of the transmission line and voltage bus dynamics. Each generator is interconnected to the network via the terminal voltage $v_\mathsf{k}$ and current $i_{s,\mathsf{k}}$ injected into the network. In the following we apply this separation to the steady-state conditions \eqref{eq:set:Sw0} and show that the steady-state conditions for the network can be used to fully characterize the synchronous steady state of the full multi-machine system. For notational convenience we define the vector of stator currents $i_s =(i_{s,1},\ldots,i_{s,n_g}) = \mathcal{I}_s i$.
The following equations describe Kirchhoff's current law at the generator terminal and load buses, as well as Kirchhoff's voltage law over the network branches:
\begin{align*}
\rho_N(i_s,v,i_T) \!=\!\!
 \begin{bmatrix}
 \!(Y_l(v)  + \omega_0 J_{v} C) v +   \!(i_s, \mathbbl{0}_{2n_l})\! + \mathcal{E} i_T   \\
 (R_T + \omega_0 J_T L_T) i_T - \mathcal{E}^{\top} v
 \end{bmatrix}\!.
\end{align*}
Based on the vector $(i_s,v,i_T) \in \mathbb{R}^{n_z}$ we define the solution set $\mathcal{N}$ of the transmission network equations
\begin{align}\label{eq:set:Nw0}
 \!\mathcal{N}_{\omega_0} \!:=\! \left\{  (i_s,v,i_T)  \in \mathbb{R}^{n_z} \left\vert\; \rho_N(i_s,v,i_T)=\mathbbl{0}_{n_N} \right.\right\}\!.\!
\end{align}
The following statement formalizes the separation that allows to recover the full system state from a solution to the steady-state equations \eqref{eq:set:Nw0} of the transmission network. In particular, given a solution to the transmission network equations $\rho_N(i_s,v,i_T)=\mathbbl{0}_{n_N}$, the remaining states and inputs such that $(x,u) \in \mathcal{S}_{\omega_0}$ can be recovered. This results in a simpler condition for steady-state operation which no longer depends on the rotor angles $\theta$ and directly links the currents $i_s$ injected by the generators to the network steady-state conditions. To this end, we define the stator impedance $Z_{s,\mathsf{k}}(\theta_\mathsf{k}) \coloneqq R_{s,\mathsf{k}} + \omega_0 j L_{s,\mathsf{k}}(\theta_\mathsf{k})$ and the voltage $\nu_\mathsf{k}(\theta_\mathsf{k}) \coloneqq v_\mathsf{k} - Z_{s,\mathsf{k}}(\theta_\mathsf{k}) i_{s,\mathsf{k}}$.

\begin{thm}{\bf(Network generator separation)}
\label{th:pr}%
Consider the sets $\mathcal{S}_{\omega_0}$ and $\mathcal{N}_{\omega_0}$ defined in \eqref{eq:set:Sw0} and \eqref{eq:set:Nw0}. If the network equations are not satisfied, i.e., $(i_s,v,i_T) \notin \mathcal{N}_{\omega_0}$, then there exist no machine states $(\theta, \omega, i)$ and no machine inputs $u$ such that $(\theta,\omega,i,v,i_T,u) \in \mathcal{S}_{\omega_0}$. Conversely, for every $(i_s,v,i_T) \in \mathcal{N}_{\omega_0}$ there exist $(\theta,\omega,i)$ and $u$ such that $(\theta,\omega,i,v,i_T,u) \in \mathcal{S}_{\omega_0}$. \\[2mm]
In particular, given  $(i_s,v,i_T) \in \mathcal{N}_{\omega_0}$ and the rotor polarization $\sigma_\mathsf{k} \in \{-1,1\}$, all corresponding $(\theta,\omega,i,u)$ such that $(\theta,\omega,i,v,i_T,u) \in \mathcal{S}_{\omega_0}$ satisfy
\begin{subequations}\label{eq:ssl}
  \begin{align}
      \omega_0 l_{sf,\mathsf{k}}  i_{f,\mathsf{k}} &= \sigma_\mathsf{k} \norm{\nu_\mathsf{k}(\theta_\mathsf{k})}, \label{eq:ssl:if1} \\
       j \mathrm{r}(\theta_\mathsf{k}) \norm{\nu_\mathsf{k}(\theta_\mathsf{k})} &= \sigma_\mathsf{k} \nu_\mathsf{k}(\theta_\mathsf{k}), \label{eq:ssl:theta1}
  \end{align}
  as well as $\omega_\mathsf{k} = \omega_0$, $i_{d,\mathsf{k}}=i_{q,\mathsf{k}}=0$, and
  \begin{align}
    v_{f,\mathsf{k}} &= r_{f,\mathsf{k}} i_{f,\mathsf{k}}, \label{eq:ssl:vf}\\
    \tau_{m,\mathsf{k}} &= D_\mathsf{k} \omega_0 + \tau_{e,\mathsf{k}}, \label{eq:ssl:tm}
  \end{align}
      where $\tau_{e,\mathsf{k}} = \tfrac{1}{2} i^{\top}_\mathsf{k} \big(L_\mathsf{k}(\theta_\mathsf{k}) \mathscrl{j}  + \mathscrl{j}^{\top} L_\mathsf{k}(\theta_\mathsf{k})\big) i_\mathsf{k}$.
\end{subequations}
\end{thm}
\begin{pf}
To prove the first statement, we note that the last two equations defining $\mathcal{S}_{\omega_0}$ are identical to the equations defining $\mathcal{N}_{\omega_0}$. It trivially follows that no $(\theta,\omega,i)$ and $u$ exist such that $(\theta,\omega,i,v,i_T,u) \in  \mathcal{S}_{\omega_0}$ if $(\mathcal{I}_s i,v,i_T) \notin \mathcal{N}_{\omega_0}$.

Next, note that $(\theta,\omega,i,v,i_T,u) \in \mathcal{S}_{\omega_0}$ requires that $\omega = \mathbbl{1}_{n_g} \omega_0$ and $(R + \omega_0 L(\theta) \mathcal{J}_g ) i = \mathcal{I}^{\top}_v v  + \mathcal{I}_f  v_f - v_{\textit{ind}}$. By using $\mathcal{J}^\top_g = - \mathcal{J}_g$ this condition simplifies to
\begin{align}
  (R + \omega_0 \mathcal{J}_g L(\theta) ) i =  \mathcal{I}_v^{\top} v  + \mathcal{I}_f  v_f.  \label{th:pr:II:3}
\end{align}
By considering the components of \eqref{th:pr:II:3}, it can be seen that this results in $i_{d,\mathsf{k}}=i_{q,\mathsf{k}}=0$, and we obtain the following condition for $i_{s,\mathsf{k}}$ and $i_{f,\mathsf{k}}$ for all 
  $\mathsf{k} \in \{1,\ldots,n_g\}$:
\begin{align*}
     \begin{bmatrix} R_{s,\mathsf{k}}+\omega_0 j L_{s,\mathsf{k}}(\theta_\mathsf{k}) &  \omega_0 j \mathrm{r}(\theta_\mathsf{k}) l_{sf,\mathsf{k}} \\ 0 & R_{f,\mathsf{k}} \end{bmatrix}  \begin{bmatrix} i_{s,\mathsf{k}} \\ i_{f,\mathsf{k}} \end{bmatrix}  = 
     \begin{bmatrix} v_\mathsf{k}  \\   v_{f,\mathsf{k}}\end{bmatrix}.
\end{align*}
This holds if and only if $v_{f,\mathsf{k}}=r_{f,\mathsf{k}} i_{f,\mathsf{k}}$ and results in:
\begin{align}
	 \omega_0 l_{sf,\mathsf{k}} i_{f,\mathsf{k}} j \mathrm{r}(\theta_\mathsf{k})  = v_\mathsf{k} - Z_{s,\mathsf{k}}(\theta_\mathsf{k}) i_{s,\mathsf{k}} = \nu_\mathsf{k}.\label{eq:vec}
\end{align}
If $\omega_0 \neq 0$ and there exist $\theta^\prime_\mathsf{k}$ such that $\nu_\mathsf{k}(\theta^\prime_\mathsf{k}) = \mathbbl{0}_2$, $i_{f,\mathsf{k}}=0$ and $\theta_\mathsf{k} = \theta^\prime_\mathsf{k}$ solve \eqref{eq:vec}. Next, we consider the non-trivial case, i.e., $\omega_0 \neq 0$ and $\nu_\mathsf{k}(\theta_\mathsf{k}) \neq \mathbbl{0}_2$ for all $\theta_\mathsf{k}$. To this end, \eqref{eq:vec} can be rewritten as follows:
\begin{align}
	  i_{f,\mathsf{k}} \eta_{f,\mathsf{k}} = \mathrm{R}(\theta_\mathsf{k})^\top \eta_{c,\mathsf{k}} +  \mathrm{R}(\theta_\mathsf{k}) \eta_{sa,\mathsf{k}},\label{eq:vec2}
\end{align}
where $\eta_{f,\mathsf{k}} \coloneqq (\omega_0 l_{sf,\mathsf{k}},0)$, $\eta_{c,\mathsf{k}} \coloneqq j^\top (v_\mathsf{k} - (R_{s,\mathsf{k}}+\omega_0j l_{s,\mathsf{k}}) i_{s,\mathsf{k}})$, and $\eta_{sa,\mathsf{k}}= \omega_0 l_{sa,\mathsf{k}} \diag(1,-1) i_{s,\mathsf{k}}$. Rewriting $\eta_{c,\mathsf{k}}$ and $\eta_{sa,\mathsf{k}}$ in polar coordinates, i.e., defining $\alpha_{c,\mathsf{k}} \in \mathbb{R}_{\geq 0}$, $\alpha_{sa,\mathsf{k}} \in \mathbb{R}_{\geq 0}$, $\delta_{c,\mathsf{k}} \in \mathbb{S}^1$, and $\delta_{sa,\mathsf{k}} \in \mathbb{S}^1$ such that $\eta_{c,\mathsf{k}} = \alpha_{c,\mathsf{k}} \mathrm{r}(\delta_{c,\mathsf{k}})$ and $\eta_{sa,\mathsf{k}} = \alpha_{sa,\mathsf{k}} \mathrm{r}(\delta_{sa,\mathsf{k}})$,  \eqref{eq:vec2} becomes $i_{f,\mathsf{k}} \eta_{f,\mathsf{k}} = \varepsilon(\theta_\mathsf{k})$ with
\begin{align}
\varepsilon(\theta_\mathsf{k}) \coloneqq \alpha_{c,\mathsf{k}} \mathrm{r}(\delta_{c,\mathsf{k}}-\theta) + \alpha_{sa,\mathsf{k}} \mathrm{r}(\delta_{sa,\mathsf{k}}+\theta).
\end{align}
If $\alpha_{c,\mathsf{k}} \!=\! \alpha_{sa,\mathsf{k}}$, the second line of $i_{f,\mathsf{k}} \eta_{f,\mathsf{k}} \!=\! \varepsilon(\theta_\mathsf{k})$ holds only if $\theta_\mathsf{k} \!=\! \theta^\star_\mathsf{k} \!=\! \tfrac{1}{2}(\delta_{c,\mathsf{k}} \!-\! \delta_{sa,\mathsf{k}} \!+\! \pi)$ and it follows that there exists $i_{f,\mathsf{k}}$ such that $i_{f,\mathsf{k}} \eta_{f,\mathsf{k}} \!=\! \varepsilon(\theta^\star_\mathsf{k})$ holds. In the following, we show that $\varepsilon(\theta_\mathsf{k})$ with $\alpha_{c,\mathsf{k}} \!\neq\! \alpha_{sa,\mathsf{k}}$ defines an ellipse with non-zero diameter, centered at the origin, and parametrized by $\theta_\mathsf{k}$. It can be verified that $\varepsilon(\theta_\mathsf{k}) + \varepsilon(\theta_\mathsf{k}+\pi) \!=\! \mathbbl{0}_2$, i.e., the ellipse is centered at the origin, and the diameter is strictly positive for all angles $\theta_\mathsf{k}$ because $\varepsilon(\theta_\mathsf{k})^{\!\top}\! \varepsilon(\theta_\mathsf{k}) \!= \alpha_{c,\mathsf{k}}^2 - 2\alpha_{c,\mathsf{k}} \alpha_{sa,\mathsf{k}} \cos(2\theta_\mathsf{k}+ \delta_{sa,\mathsf{k}} - \delta_{c,\mathsf{k}}) + \alpha_{sa,\mathsf{k}}^2 \geq (\alpha_{c,\mathsf{k}} - \alpha_{sa,\mathsf{k}})^2$ holds. It follows that there exist two $\theta_\mathsf{k}$ such that the second component of $\varepsilon(\theta_\mathsf{k})$, i.e., of the right hand side of \eqref{eq:vec2}, is equal to zero. Moreover, if $\omega_0 \neq 0$ there exist pairs $(i_{f,\mathsf{k}},\theta_\mathsf{k}) \in \mathbb{R}_{\geq 0} \times \mathbb{S}^1$ and $(i_{f,\mathsf{k}},\theta_\mathsf{k}) \in \mathbb{R}_{\leq 0} \times \mathbb{S}^1$ such that $i_{f,\mathsf{k}} \eta_{f,\mathsf{k}} = \varepsilon(\theta_\mathsf{k})$ is satisfied. Considering \eqref{eq:vec}, any such pair satisfies \eqref{eq:ssl:if1} and \eqref{eq:ssl:theta1}.

Next, $\omega_0=0$ results in the steady-state condition $\ddt i_{s,\mathsf{k}} = \mathbbl{0}_2$ and $v_{\textit{ind},\mathsf{k}}=\mathbbl{0}_5$. Considering \eqref{eq:gendyn:elec} this results in $\nu_\mathsf{k} = v_\mathsf{k} - R_{s,\mathsf{k}} i_{s,\mathsf{k}} = \mathbbl{0}_2$. In other words, $\omega_0=0$ implies $\nu_\mathsf{k} = \mathbbl{0}_2$ in steady state for all $\theta_\mathsf{k}$, and \eqref{eq:vec}, \eqref{eq:ssl:if1}, and \eqref{eq:ssl:theta1}, are satisfied for any $\theta_\mathsf{k}$ and any $i_{f,\mathsf{k}}$.

Finally, note that $(\theta,\omega,i,v,i_T,u) \in \mathcal{S}_{\omega_0}$ requires $\tau_m = D \omega + \tau_e$. Using the currents $i_{s,\mathsf{k}}$, $i_{f,\mathsf{k}}$, and the angle $\theta_\mathsf{k}$, the electrical torque $\tau_{e,\mathsf{k}}$ can be explicitly recovered using \eqref{orig:torque} resulting in \eqref{eq:ssl:tm}. Moreover, the steady-state excitation voltage $v_{f,\mathsf{k}}$ is given by \eqref{eq:ssl:vf}. It follows that for every $(i_s,v,i_T) \in \mathcal{N}_{\omega_0}$ there exist $(\theta,\omega,i)$ and $u$, such that $(\theta,\omega,i,v,i_T,u) \in \mathcal{S}_{\omega_0}$. \hfill $\square$
\end{pf}
Broadly speaking, Theorem \ref{th:pr} shows that the conditions for operation of the power grid \eqref{eq:nlmodel} at the synchronous steady state can be equivalently expressed
in terms of the transmission network equations and the currents injected into the network by the generators. In particular, for every solution of the network equations $\rho_N(i_s,v,i_T)=\mathbbl{0}_{n_N}$ one directly obtains the function $\nu_\mathsf{k}(\theta_\mathsf{k})$ for all $\mathsf{k} \in \{1,\ldots,n_g\}$. Moreover, given $\nu_\mathsf{k}(\theta_\mathsf{k})$, field winding currents $i_{f,\mathsf{k}}$, and angles $\theta_\mathsf{k}$ which satisfy the steady-state conditions can be computed using \eqref{eq:vec}. Specifically, if $\omega_0 \neq 0$ and $l_{sa,\mathsf{k}} = 0$ for all $\mathsf{k} \in \{1,\ldots,n_g\}$, $\nu_\mathsf{k}$ is independent of $\theta_\mathsf{k}$ and one obtains $ i_{f,\mathsf{k}} = (\omega_0 l_{sf,\mathsf{k}})^{-1} \sigma_\mathsf{k} \norm{\nu_\mathsf{k}}$ with $\sigma_\mathsf{k} \in \{-1,1\}$ and a corresponding $\theta_\mathsf{k}$ can be directly computed from \eqref{eq:ssl:theta1}. Together with $i_{d,\mathsf{k}}=i_{q,\mathsf{k}}=0$ this recovers a steady state of the full system from a solution to the network equations. Moreover, the corresponding steady-state inputs $u$ can be computed as follows. Given a steady state of the generators, the electrical torque $\tau_{e,\mathsf{k}}$ can be explicitly recovered using \eqref{orig:torque} and the mechanical torque $\tau_{m,\mathsf{k}}$ can be computed using \eqref{eq:ssl:tm}. Finally, the excitation voltage $v_{f,\mathsf{k}}$ is given by \eqref{eq:ssl:vf}. We emphasize that there exist \emph{two} steady-state solutions for the angle of \emph{each} generator (cf. \cite{BSO+16} for a similar result for a single generator connected to an infinite bus).
\begin{rem}{\bf(Trivial cases)}
If $\omega_0 \neq 0$ and $\nu_\mathsf{k}(\theta_\mathsf{k}) = \mathbbl{0}_2$ for the steady-state angle satisfying \eqref{eq:ssl}, the conditions \eqref{eq:ssl} hold for $i_{f,\mathsf{k}}=0$, $v_{f,\mathsf{k}}=0$ and any angle $\theta_\mathsf{k}$. The mechanical torque $\tau_{m,\mathsf{k}}=D_\mathsf{k} \omega_0$ only compensates for the losses in the generator. In contrast, $\omega_0 = 0$ implies $\nu_\mathsf{k}=\mathbbl{0}_2$ for any any $\theta_\mathsf{k}$ and the steady-state conditions hold for any $i_{f,\mathsf{k}}$ and any angle $\theta_\mathsf{k}$. In this case, $\tau_{m,\mathsf{k}}=\tau_{e,\mathsf{k}}$ holds the rotor in place (i.e., $\omega_\mathsf{k}=0$). Clearly these cases do not define sensible operating points. \hfill $\square$
\end{rem}
In the next section, we establish a connection between the results of Theorem \ref{th:pr} and nodal current balance (or power flow) equations typically used to determine the desired operating point of power systems. 

\subsection{Nodal Balance Equations}
In practice, the desired operating point of power systems is determined based on nodal current balance (or power flow) equations and found via, e.g., load flow analysis and generation dispatch optimization. In the following, we show that the nodal current balance equations fully characterize the steady-state behavior of the overall power system, i.e. for every solution of the nodal balance (or power flow) equations there exists a corresponding synchronous steady-state behavior.

The set $\mathcal{P}_{\omega_0}$ describes the nodal current balance equations \cite{K94}:
\begin{align}\label{eq:nodalbalance}
 \mathcal{P}_{\omega_0} := \left\{ (i_s,v) \in \mathbb{R}^{n_p} \left\vert \;
 (-i_s, \mathbbl{0}_{2n_l}) = Y_N v \right.\right\}.
\end{align}
Here $Y_v = Y_l(v) + \omega_0 J_{v} C$ is the matrix of shunt load admittances, 
$Z_T = R_T + \omega_0 J_T L_T$ is the matrix of branch impedances, $Y_N = Y_v\! + \mathcal{E} Z^{-1}_T \mathcal{E}^{\top}$ 
is the network admittance matrix, and $Z_T$ is invertible for all $\omega_0 \in \mathbb{R}_{\geq 0}$. 
\begin{thm}{\bf(Nodal current balance equations)}\label{th:nodalbalance}
Consider the sets $\mathcal{S}_{\omega_0}$ and $\mathcal{P}_{\omega_0}$ defined in \eqref{eq:set:Nw0} and \eqref{eq:nodalbalance}. If the nodal balance equations are not satisfied, i.e., $(i_s,v) \notin \mathcal{P}_{\omega_0}$, then there exist no states $(\theta, \omega, i, i_T)$ and no machine inputs $u$ such that $(\theta,\omega,i,v,i_T,u) \in \mathcal{S}_{\omega_0}$. Conversely, for every $(i_s,v) \in \mathcal{P}_{\omega_0}$ there exist $(\theta,\omega,i,i_T)$ and $u$ such that $(\theta,\omega,i,v,i_T,u) \in \mathcal{S}_{\omega_0}$.
\end{thm}
\begin{pf}
Consider $i_T = Z^{-1}_T \mathcal{E}^{\top} v$. It follows that $(i_s,v,Z^{-1}_T \mathcal{E}^{\top} v) \notin \mathcal{N}_{\omega_0}$ for all $(i_s,v) \notin \mathcal{P}_{\omega_0}$. Moreover, by definition of $\mathcal{N}_{\omega_0}$ it follows that there cannot exist $i_s$ such that $(i_s,v,i_T)\in\mathcal{N}_{\omega_0}$ for $i_T \neq Z^{-1}_T \mathcal{E}^{\top} v$. In other words, for $(i_s,v) \notin \mathcal{P}_{\omega_0}$, no $i_T$ exists such that $(i_s,v,i_T)\in\mathcal{N}_{\omega_0}$. Considering Theorem \ref{th:pr}, it directly follows that there exist $(\theta, \omega, i)$ and $u$ such that $(\theta,\omega,i,v,i_T,u) \in \mathcal{S}_{\omega_0}$ only if $(i_s,v) \in \mathcal{P}_{\omega_0}$. Next, for every $(i_s,v) \in \mathcal{P}_{\omega_0}$ it holds that $(i_s,v,Z^{-1}_T \mathcal{E}^{\top} v) \in \mathcal{N}_{\omega_0}$. Considering Theorem \ref{th:pr}, it follows that for every $(i_s,v) \in \mathcal{P}_{\omega_0}$ there exist $(\theta, \omega, i)$ and $u$ such that $(\theta,\omega,i,v,i_T,u) \in \mathcal{S}_{\omega_0}$ holds. \hfill $\square$
\end{pf}
Theorem \ref{th:nodalbalance} shows that for every non-trivial solution (i.e., with non-zero currents and voltages) to the nodal current balance equations there exists a corresponding
non-trivial steady-state behavior of the power system.
\begin{rem}{\bf(Power flow equations in polar coordinates)}
 The standard power flow equations (studied in, e.g., \cite{DBO+09}) can be recovered from $\mathcal{P}_{\omega_0}$ by multiplying the nodal balance equations \eqref{eq:nodalbalance} by $v$ from the left and rewriting the resulting equations  in complex phasor notation. \hfill $\square$
\end{rem}

\section{Main Result and Discussion}\label{sec:disc}
\begin{thm}{\bf(Steady-state operation)}\label{cor:constinp}
 Consider the power system dynamics \eqref{eq:nlmodel}, a solution $(i_s,v) \in \mathcal{P}_{\omega_0}$ to the nodal current balance equations, and $i_{T} = Z^{-1}_T \mathcal{E}^{\top} v$. The power system operates in steady state if and only if:\\[-7mm]
 \begin{enumerate}
  \item the machine state $(\theta,\omega,i)$, the input $u$, and $(i_s,v)$ satisfy \eqref{eq:ssl},
  \item the input $u$ is constant,
  \item the load model can be written as $i_{l,\mathsf{k}} = Y_l(\norm{v_\mathsf{k}}) v_\mathsf{k}$. 
 \end{enumerate}
\end{thm}
Theorem \ref{cor:constinp} combines necessary and sufficient conditions to ensure steady-state operation of the detailed first-principles power system model proposed in \cite{FZO+13}. It shows that several steady-state conditions which are commonly assumed to be sufficient in the analysis of power system are, in fact, necessary and sufficient for synchronous and balanced steady-state operation. In other words, any power system control strategy and any load model necessarily need to satisfy the conditions given in Theorem \ref{cor:constinp}. Conventional dispatch optimization and generator control strategies such as speed droop control, automatic voltage regulators, all implicitly meet these specifications in steady state. Moreover, the steady-state specifications for the generators justify assumptions, such as constant excitation current, which are often made in the stability analysis of multi-machine networks \citep{CT14} even when the power system is not in steady state.

Our analysis and rigorous definition of synchronous balanced steady-states via the set $\mathcal{S}$ can be seen as starting point for multi-machine power system stability analysis and control design which avoids the inherent difficulties of using local rotating coordinate frames (i.e., local $dq$ frames) for each device. While such rotating coordinate frames are convenient for a single generator they severely complicate the stability analysis for multiple generators \citep{CT:16}. In contrast, our definition of a synchronous steady state does not require rotating coordinates frames. Moreover, we fully characterize the class of static load models which are compatible with synchronous and balanced steady-state operation and thereby narrow down the class of load models to be considered in stability analysis. Finally, the approach used in this work is directly applicable to characterize the steady-state and steady-state control inputs of low-inertia power systems with renewable generation interfaced by power electronics \citep{GD17}.

\section{Conclusion}\label{sec:concl}
In this paper, we provided results on the steady-state behavior of a nonlinear multi-machine three-phase power system model including nonlinear generator dynamics, a dynamic model of the transmission network, and static nonlinear loads. The steady-state behavior considered in this work is defined by balanced and sinusoidal three-phase AC signals of the same synchronous frequency.
In the literature on power systems it is often assumed a priori that the power system admits such a steady-state behavior if the field current and mechanical torque input are constant, the nodal current balance (or power flow) equations can be solved, and specific load models are used. We show that all of these conditions can be constructively obtained from first-principle and are, in fact, necessary and sufficient for the power system to admit synchronous and balanced steady-state behaviors. Extending the results to include power converters  gives rise to non-trivial internal models for power converter control and is the focus of ongoing work.

\bibliographystyle{elsarticle-harv} 
\bibliography{steadystate}   

\end{document}